# Distributed Optimal Vehicle Grid Integration Strategy with User Behavior Prediction

Yingqi Xiong, *Student Member, IEEE*, Chi-cheng Chu
and Rajit Gadh
Department of Mechanical Engineering
University of California, Los Angeles
Los Angeles, USA
yxb936@ucla.edu

Bin Wang
Grid Integration Group
Lawrence Berkeley National Laboratory
Berkeley, USA
wangbin@lbl.gov

*Abstract*—**With the increasing of electric vehicle (EV) adoption in recent years, the impact of EV charging activities to the power grid becomes more and more significant. In this article, an optimal scheduling algorithm which combines smart EV charging and V2G gird service is developed to integrate EVs into power grid as distributed energy resources, with improved system cost performance. Specifically, an optimization problem is formulated and solved at each EV charging station according to control signal from aggregated control center and user charging behavior prediction by mean estimation and linear regression. The control center collects distributed optimization results and updates the control signal, periodically. The iteration continues until it converges to optimal scheduling. Experimental result shows this algorithm helps fill the valley and shave the peak in electric load profiles within a microgrid, while the energy demand of individual driver can be satisfied.**

*Index Terms*—**Distributed Optimization, Electric Vehicle Charging, Vehicle to Grid, User Behavior Prediction.**

## I. Introduction

Global warming concerns make Electric Vehicle (EV) and Plug-in Hybrid Electric Vehicle (PHEV) more and more popular in recent years. It is estimated that by 2020, EVs will comprise about 4 percent of the 100 million light-duty vehicles [1]. Four million EVs in the US will produce 40 GW of potentially available power if each EV is grid-integrated with an average of 10 kW of available capacity. Currently there are lack of policies to regulate EV charging. All charging activities are spontaneous without coordination and considerations of other grid circumstances. Un-coordinated EV charging behaviors will impact the power grid and degrade power quality when the number of charging sessions in local distribution network reaches a certain level [2], resulting in voltage fluctuation, power loss and blackout, etc. Smart charging can solve these issues by optimizing the EV charging schedule and energy allocation [3]-[5]. Significant studies have been conducted to incorporate bi-directional power flow and grid services support from EVs, upgrading traditional dumb consumptions to controllable distributed energy resources (DERs) [6]. If one can combine EV charging and V2G service into a smart bi-directional charging strategy, EVs will become valuable assets to the power grid. However, distributed EV batteries must be aggregated in a large scale so that they can participate in the wholesale market to provide grid services. Such aggregation would require huge computational resources if it's performed by traditional centralized control scheme. Smart bi-directional charging scheduling depends highly on EV driver behaviors. To facilitate the integration of EVs as DERs into power grid, the following challenges should be addressed: 1) A bi-directional charging algorithm incorporates both smart charging and V2G service; 2) A decentralize control strategy which lowers computational cost with large-scale integration; 3) EV user behavior modeling and prediction.

Previous research has been done to patricianly solve EV charging coordination problems. But none of them provides a comprehensive solution which addresses all the three aforementioned challenges based on real-world implementation. Gan et al. [5] use optimal control to optimally allocate EV charging time and energy, but the algorithm requires users to provide charging schedule, which may degrade customers' satisfaction by requiring inputs frequently. Wang et al. [7] provide two algorithms to aggregate EV charging with consideration of random user behavior model and renewable generation in the scheduling. Authors in [8] have developed smart charging strategy according to time-of-use (TOU) price from day-ahead predictions. These smart charging strategies can flatten the grid load profile by valley-filling but cannot do much during peak load. V2G can be performed to relieve the stress on the grid [9]-[11]. Methods for increasing smart charging and V2G scalability and interoperability are developed in [10], but large scale implementation strategy based on real-world data is rarely investigated.

In this paper we proposed a smart bi-directional charging scheduling algorithm that incorporates EV smart charging and V2G service for large scale implementation. An aggregated control center sends control signal to all EV charging stations in its network. Each distributed charging station performs an optimization locally according to this control signal to update

This work has been sponsored in part by grants from the California Energy Commission (CEC) entitled "Demonstration of PEV Smart Charging and Storage Supporting Grid Operational Needs". Sponsor Award number: EPC-14-056

EV bi-directional charging strategy and reports it to control center. Subsequently, control center collects strategies from networked charging stations and then updates its original control signal. The iteration ends when algorithm converges to an optimal charging strategy for all plug-in EVs in the network. Computational cost is shared by all charging stations which greatly lowers computing burden and time in control center. In each iteration, the information update is asynchronous, there is no need for all charging stations update their strategy. Also there is no strict rule for control center to update control signal every iteration. These features are critical for a successful large scale implementation where prompt updating is impractical. User charging behavior and preference prediction is performed by mean estimation and linear regression using real-world EV usage data. Predicted EV user charging schedule and energy demand serve as input parameters for the proposed algorithm. The contributions of this paper can be summarized as: 1) A distributed optimal bi-directional charging control algorithm implementation with real-world data; 2) Integrate EVs as DERs into power grid to flatten load curves, correct voltage deviation, conduct voltage regulation, etc.; 3) A driver behavioral model based on mean estimation and linear regression is developed.

The rest of this paper is organized as follows: Section II discusses the smart charging network infrastructure for proposed algorithm implementation and EV user behavior prediction. Section III covers the details of proposed algorithm. Section IV shows implementation result and evaluates the performance.

## II. SYSTEM OVERVIEW

### A. Smart Charging Infrastrcuture

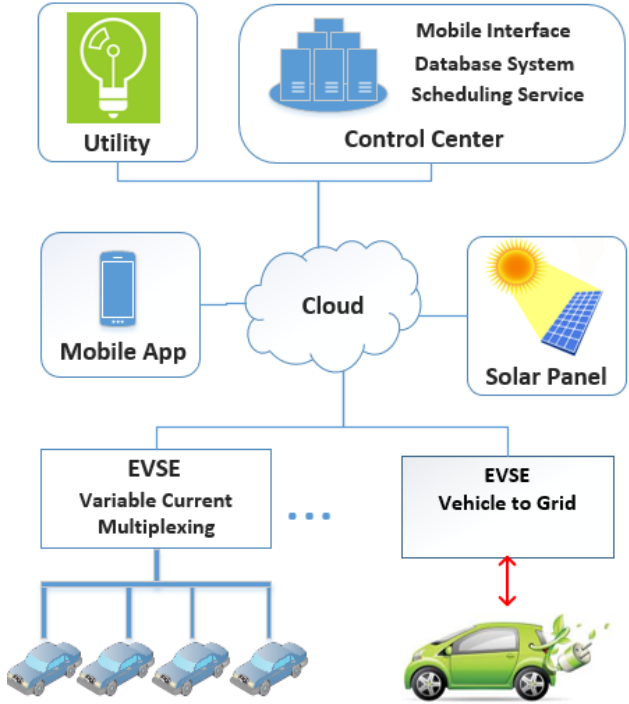

Figure 1. Smart charging system architecture

The smart charging infrastructure, shown in Fig. 1, is the testbed on UCLA campus to implement proposed bi-directional charging control algorithm. Developed by pioneer researchers in UCLA Smart Grid Energy Research Center (SMERC), this infrastructure has already supported 17 smart electric vehicle supply equipment (EVSEs) around Los Angeles area. Smart EVSEs are the foundations of this infrastructure. All of them are connected with control center through 3G network. Each smart EVSE has an embedded charging control system coordinating local charging session with variable current charging and power sharing feature [12]. Smart charging algorithm [13]-[15] considering the energy sharing strategies, uncertainties of user behaviors and renewable solar generations has been successfully implemented in this infrastructure with real-world EV users. EVSEs with V2G feature have been added to this smart charging network.

### B. EV User Behavior Prediction

The smart charging infrastructure provides mobile application for its users to specify their energy demand and charging preference. Control center records user charging data and stores them in its database. In this paper, charging data in the past 3 years are extracted from database for user behavior analysis.

Historical charging data of randomly selected 4 EV users are plot in Fig. 2 to show the relationship between charging start time and end time.

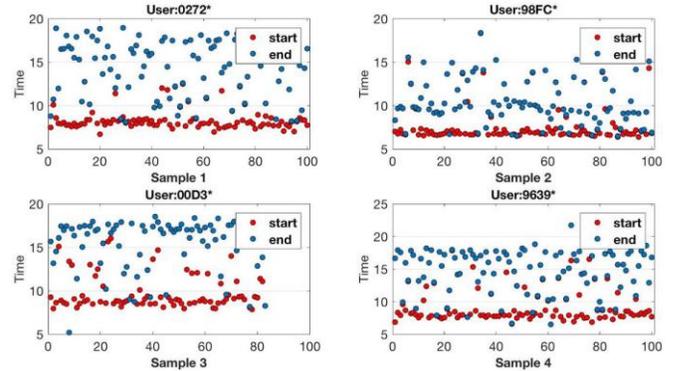

Figure 2. EV user charging time preference

From Fig. 2 it can be seen that the charging start-time and end-time of these 4 EV drivers follow some patterns. Most of the charging activities start and stop at particular time in the day, which may represent the EV user weekdays routine. To predict the charging behavior of these EV users, mean estimation is used. The charging start time and end time can be obtained using follow equations:

$$t_{start}^{pred} = \frac{1}{M}\sum_{i=1}^{M} t_{start}^{i} \qquad (1)$$

$$t_{end}^{pred} = \frac{1}{M}\sum_{i=1}^{M} t_{end}^{i} \qquad (2)$$

where $t_{start}^{pred}$ is the predicted charging start time using mean method, and $t_{end}^{pred}$ is the end time. M is the number of samples used in prediction analysis.

EV charging duration $t_{stay}$ can be obtained by taking the difference between start time and end time. Fig. 3 shows the relationship between charging duration and energy consumption for one of the users.

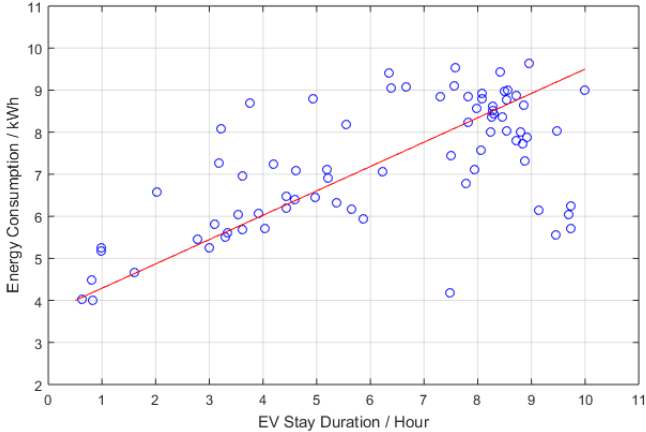

Figure 3. EV charging power consumption versus duration

It can be seen from Fig. 3 the relationship between charging duration and energy consumption is almost linear. Additional features, such as last trip energy consumption, driver' home distance from campus, etc, will be utilized in the regression model in the future. By using linear regression, the prediction model for EV charging consumption can be obtained. Let $X = [t_{stay}^1, t_{stay}^2, \ldots t_{stay}^M]$ denote the pool of sample charging duration and $y = [E^1, E^2, \ldots E^M]$ contains corresponding energy consumption, linear regression parameter $\theta$ can be calculated by:

$$\theta = [X^T X]^{-1} X^T y \quad (3)$$

The predicted energy consumption can then be calculated by:

$$E^{pred} = \theta(t_{end}^{pred} - t_{start}^{pred}) \quad (4)$$

Actually most of the charging behaviors from different EV drivers in our database share similar patterns with EV drivers in Figure 3. They have routine charging start time and end time with predictable energy consumption values. By using mean estimation and linear regression method, we can predict EV user charging behaviors for the proposed bi-directional charging control algorithm.

### III. CHARGING CONTROL ALGORITHM

In this section, details of the optimal bi-directional charging control algorithm are discussed, including its constraints, objective and convergence.

#### A. Problem Formulation

It is assumed that, the utility or the microgrid operator can always provide day-ahead prediction of the load profile in its service area or within the microgrid. Such prediction is defined as baseload in the scope of this paper. Control center of the smart charging infrastructure can retrieve baseload data from service provider, i.e. utility or microgrid operator as input for its optimal bi-directional charging control algorithm. Control center makes day-ahead prediction of user charging behaviors to obtain predicted charging start time, end time, energy demand for each EV users within its control network. Predicted data and control signal are then distributed to each networked EVSE for bi-directional charging profile optimization. Updated individual EV charging profiles are returned to control center for control signal updating in a new iteration. All charging stations will interact with the centralized control center and cooperatively solve the optimization problem in a distributed fashion. The goal of the algorithm is to schedule the time and allocate the power flow for charging and V2G such that the load profile in the utility service area or within the microgrid are flattened and stabilized while all EV user demands are satisfied.

#### B. Constraints and Objective

The day-ahead predicted baseload is denoted as $B(t)$ with $t \in T$. $T = [1, 2, \ldots T]$ is the time slot set where the algorithm will be performed and $t$ is the specific time slot. We define there are $N$ EVs in total within network. Each EV has a charging/discharging rate $p_n(t)$ with $n \in N = [1, 2, \ldots N]$ at time slot $t$. The aim of proposed algorithm is to flatten the total load profile:

$$L(t) = \left(B(t) + \sum_{n=1}^{N} p_n(t)\right)^2 \quad (5)$$

The charging rate and V2G capacity of EVs made by different manufacturers are varied. We use $\overline{p_n}$ and $\overline{d_n}$ to denote maximum charging rate and maximum V2G capacity of EV $n$, respectively. V2G discharging rate $\overline{d_n}$ is a negative value since the power flow is reversed. Thus we have:

$$\overline{d_n} \leq p_n(t) \leq \overline{p_n} \quad (6)$$

In our proposed control scheme, $\overline{p_n}$ and $\overline{d_n}$ are regarded as upper bound and lower bound for the controllable charging rate. During the time between EV plug in EVSE, $t_{start}^{pred}$ and EV leave, $t_{end}^{pred}$, charging rate can be set continuously within upper and lower bound, and it will be 0 otherwise:

$$\overline{p_n} = \begin{cases} p_n^{\max}, & t \in [t_{start}^{pred}, t_{end}^{pred}] \\ 0, & otherwise \end{cases} \quad (7)$$

$$\overline{d_n} = \begin{cases} d_n^{\max}, & t \in [t_{start}^{pred}, t_{end}^{pred}] \\ 0, & otherwise \end{cases} \quad (8)$$

Equation (7) and (8) are corresponding to the situation that when an EV is plugged in, its charging rate can range from maximum V2G discharging rate to maximum charging rate.

When EV leaves the EVSE, all rates are zero which means there is no existing charging or discharging.

Energy consumptions of each EV are predicted using the aforementioned user behavior model. The predicted value $E_n^{pred}$ represents energy demand of the particular EV user during his or her charging sessions. The summation of the product of EV charging rate and time interval equals to the energy demand:

$$\sum_T p_n(t)\Delta T = E_n^{pred} \quad (9)$$

where $\Delta T$ is constant since each time slot is a constant time value.

The control signal from control center is the derivative of total load profile we intent to flatten multiplied by a tuning parameter [5]:

$$c^i(t) = \frac{1}{\lambda N}\left(B(t) + \sum_{n=1}^{N} p_n(t)\right) \quad (10)$$

where $i$ represents the iteration number in the algorithm and $\lambda$ denotes the control parameter. The algorithm converges when the difference between control signals in consecutive iterations smaller than a convergence criteria, i.e. $\|c^{i+1} - c^i\| \leq \varepsilon$, where $\varepsilon$ is a positive real number small enough to be the convergence criteria. It is not necessary to update control signal in every iteration for an optimal convergence.

Subsequently, each EVSE carry out an optimization of its plug-in EV locally, responding to the control signal. The objective is as following:

$$\sum_{t=1}^{T} c^i(t) p_n^{i+1}(t) + \frac{1}{2}\|p_n^{i+1}(t) - p_n^i(t)\|^2 \quad (11)$$

By minimizing the objective equation, EVSE will obtain an updated EV charging profile $p_n^{i+1}(t)$ and report it to the control center. Again, it is not necessary to report updated charging profile to the control center in every iteration.

### C. Optimal Distributed Bi-directional Charging Algorithm

The distributed algorithm to optimize EV charging energy allocation and scheduling are presented as following:

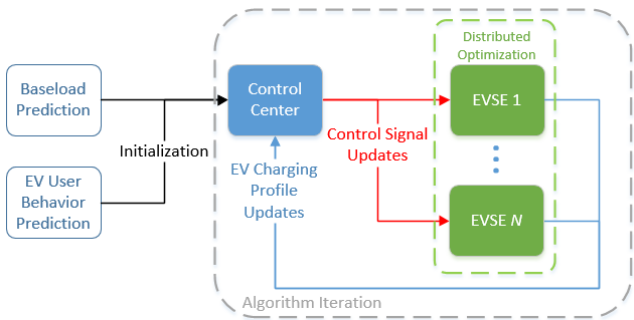

Figure 4. Schematic of proposed charging control algorithm

**Algorithm : Optimal Distributed Bi-directional Charging**

Obtain Day-ahead prediction of baseload data
Perform Day-ahead prediction of EV user behavior and energy demand
Initialize all EV charging profiles $p_n^0(t) = 0$, $n = 1,2,...N$
Pick control parameter $\lambda$
Initialize control signal $c^0(t) = \frac{1}{\lambda N}\left(B(t) + \sum_N p_n^0(t)\right)$
**While** $\|c^{i+1} - c^i\| \geq \varepsilon$ :
    **For** each $n$ EVSE in the network, $n = 1,2,...N$
        Minimize (11)
        Subject to (6), (7), (8), (9)
    **End**
    **If** $i$ mod $u == 0$
        Update control signal
    **End**
    **If** $i$ mod $v == 0$
        Update charging profile;
    **End**

In the above algorithm, Constant values u and v are used to emulate delay updating of control signal and charging profile in real-world large scale implementation. When the algorithm converges, the total load profile defined in (5) will be effectively flattened. Fig. 4 shows the schematic of proposed algorithm. Each EVSE performs optimization in its own embedded system, which greatly lowers the computational burden of the control center.

### IV. RESULTS AND DISCUSSION

This section covers the implementation of proposed algorithm on real-world EV users in UCLA microgrid and its performance evaluation.

#### A. Experiment Setup

UCLA Engineering IV building load profile on September 20, 2016 from 7 am to 7 pm is chosen as the baseload for proposed algorithm implementation. The 12 hour time span is divided into 60 time slots. The historical data in database of UCLA SMERC smart EV charging system, including 30 EV drivers on UCLA campus, are extracted for user behavior modeling and prediction. There are power usage peak and valley in the load profile corresponding to the operation of some heavy energy consuming devices in the building. The fleet of 30 EVs will provide grid service to flatten the load under the control of proposed algorithm. Such grid service can be possibly scaled up to support a microgrid or a utility service area with enough EV participation.

#### B. Performance Evaluation

Before implementation of the proposed algorithm, EV user behavior predictions are made based on the random selected historical charging record data from 30 EV users in the past 3 years. Charging start time, end time and energy demand are predicted for the day of September 20, 2016 and then incorporated in the proposed algorithm. The algorithm converges to an optimal bi-directional charging strategy which effectively flattens the original baseload by peak shaving and valley filling, as shown in Fig. 5.

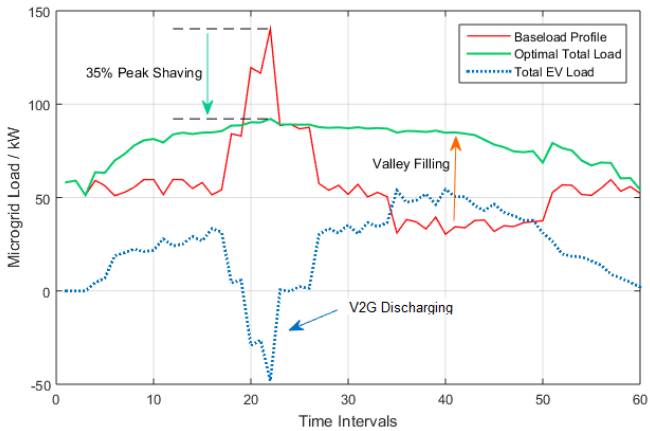

Figure 5. Base load profile flattened by proposed algorithm

In Fig. 5 it can be seen that the baseload has power consumption peak from time slot No. 17 to time slot No. 27, and power usage valley from time slot No. 33 to No. 50. The optimal distributed bi-directional charing algorithm precisely tunes charing rate of each EV in the network, changing from high speed charing to discharging according to the trend of baseload profile. Almost all EVs are performing V2G discharing at baseload peak time, making the peak power consumption drops 35%, from 140 kW down to 90 kW. The Optimal total load curve, which is the combination of baseload and EV load profile, has been perfectly flattened. The optimal distributed bi-directional charging algorithm shows its capability to integarte EVs into power grid as DERs, providing various grid service to benefits the power grid.

The values of convergence criteria ε during algorithm iteration are plot in Fig. 6. It should be notice zero is not an indication of convergence since control signal and EV charing profiles are not updated in every iterations. The algorithm converges at around 32 iterations since ε becomes enough small. This asynchronous feature of proposed algorithm shows its potential for large scale implememtation where information transmisson and updates are usually delayed.

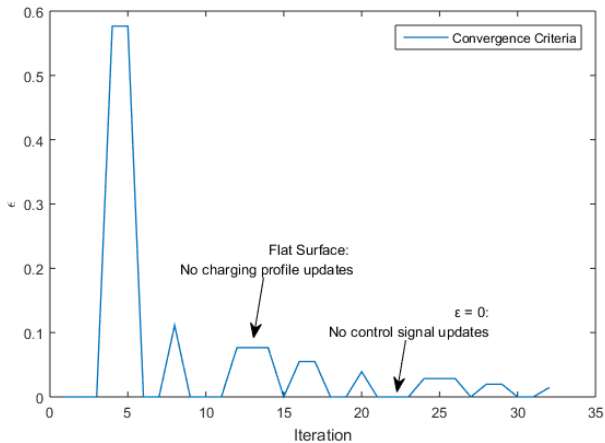

Figure 6. Convergence criteria in iterations

## V. CONCLUSION

In this paper, an optimal distributed bi-directional charging control algorithm is developed. It is implemented with real-world load profile and EV charging data. The performance of the algorithm is evaluated. Result shows its capability of integrating EVs into power grid as DERs to perform grid service, flatten the load profile and stabilize the power grid. The algorithm has distributed and asynchronous features which greatly low the computational cost and ideal for large scale implementation.